\documentclass[12pt]{article}

\usepackage{amsmath}
\usepackage{amsthm}
\usepackage{amsfonts}
\usepackage{amssymb}
\usepackage[all]{xy}
\usepackage{amsmath,amscd}

\newcommand{\R}{\mathbb{R}}

\newcommand{\PP}{\mathbb{P}}
\newcommand{\RR}{\mathbb{R^{+}}}
\newcommand{\Bl}{\operatorname{Bl}}

\newcommand{\NE}{\operatorname{NE}}

\newcommand{\Exc}{\operatorname{Exc}}

\newcommand{\OO}{{\cal O}}

\newcommand{\wtilde}{\widetilde}

\newtheorem{theo}{Theorem}

\newtheorem{lem}{Lemma}

\newtheorem{claim}{Claim}

\title
{Fano manifolds obtained by blowing up along curves with maximal Picard number}
\author{Toru Tsukioka}

\begin{document}

\maketitle

\begin{abstract}
The Picard number of a Fano manifold $X$ obtained by blowing up a curve in a smooth projective variety 
is known to be at most 5, in any dimension greater than or equal to 4. 
We show that the Picard number attains to the maximal
if and only if $X$ is the blow-up of the projective space whose center consists of two points, the strict transform of 
the line joining them 
and a linear space or a quadric of codimension 2. 
This result is obtained as a consequence of a classification of  special types of Fano manifolds.
\end{abstract}

\section{Introduction}

Let $X$ be a Fano manifold obtained by blowing up along 
a curve, i.e. there exists a pair $(Y,C)$ of 
a smooth projective variety $Y$ and 
a smooth connected curve $C\subset Y$   
such that the anticanonical divisor $-K_{X}$ is ample.
Using a recent result on Minimal Model Program due to \cite{BCHM}, 
C. Casagrande shows that such a Fano manifold has Picard number at most 5 (see \cite{Cas} Theorem 4.2 for a more 
general statement, and see \cite{S} for the toric case). 

The purpose of the present paper is to classify the maximal case:  
\begin{theo}\label{main}
Let $Y$ be a smooth projective variety of dimension $n\geq 4$ defined over the field of complex numbers, $C$ a smooth curve on $Y$, and 
$X$ the blow-up of $Y$ along $C$.  
Assume that $X$ is a Fano manifold and has Picard number 5. Then, the pair $(Y,C)$ is exactly one of the following:
\begin{enumerate}
\item $Y$ is the blow-up of $\PP^{n}$ whose center is the 
union of two points $p,q$ and 
$\PP^{n-2}$ disjoint from $\overline{pq}$, and $C$ is the strict transform of $\overline{pq}$,
\item $Y$ is the blow-up of $\PP^{n}$ whose center is the 
union of two points $p,q$ and a smooth quadric 
$Q_{n-2}$ disjoint from $\overline{pq}$, and $C$ is the strict transform of $\overline{pq}$. 
\end{enumerate}
\end{theo}
\noindent{\it Remark.} We denote by $\overline{pq}$ 
the line passing through $p$ and $q$ in $\PP^{n}$.

\

According to Casagrande's result (see \cite{Cas} Theorem 4.2 (ii)), if the assumption of Theorem \ref{main} 
is satisfied, then there exists another structure of 
blow-up $\varphi: X\to Z$ with the following properties: 
\begin{itemize}
\item $Z$ is a smooth projective variety, and 
the center of the blow-up $\varphi$ is a smooth subvariety 
of codimension 2 
\item $E\cdot f>0$, where $E$ is the exceptional 
divisor of the blow-up $\pi:X\to Y$ and $f$ is a non trivial fiber of $\varphi$
\item $F\cdot e=0$, where $F$ is the exceptional 
divisor of $\varphi$ and $e$ is a line in a fiber of 
the $\PP^{n-2}$-bundle $\pi |_{E}:E\to C$
\end{itemize} 
Hence, our Theorem \ref{main} is a consequence of the following classification result (in which only two examples (8) and (9) have Picard number 5):  
\begin{theo}\label{Fe=0} Let $Y$ be a complex manifold of dimension $n\geq 4$. Assume that 
there exists a smooth curve $C\subset Y$ such that the blow-up $X$ of $Y$ along $C$ is 
a Fano manifold. Assume moreover that 
there exists a smooth projective variety 
$Z$ and a smooth subvariety $W\subset Z$ 
of codimension 2 such that the blow-up of $Z$ along $W$ is 
isomorphic to $X$. 
Let $E$ (resp. $F$) be the exceptional divisor of the blow-up $\pi:X\to Y$ (resp. $\varphi:X\to Z$).  
Let $e$ (resp. $f$) be a line in a fiber of 
the $\PP^{n-2}$-bundle $\pi |_{E}:E\to C$ 
(resp. a fiber of the $\PP^{1}$-bundle 
$\varphi |_{F}:F\to W$). 
If $E\cdot f>0$ and $F\cdot e=0$, 
then we have exactly one of the following: 

\begin{enumerate}
\item  $Y$ is the blow-up of $\PP^{n}$ at a point $p$
and  $C$ is the strict transform of the line 
passing through $p$,
\item $Y$ is the blow-up of $Q_{n}$ at a point $p$
and  $C$ is the strict transform of a line passing through $p$,
\item  $Y$ is the blow-up of $Q_{n}$  
at a point $p$ and $C$ is the strict transform of a conic passing through $p$,
\item $Y$ is the blow-up of $\PP^{n}$ 
whose center is the union of a point $p$ and 
a linear subspace $P\simeq \PP^{n-2}$ 
not containing $p$, 
and $C$ is the strict transform of 
a line passing through $p$ and disjoint from $P$,
\item $Y$ is the blow-up of $\PP^{n}$ 
whose center is the union of a smooth 
quadric $Q\simeq Q_{n-2}$ and 
a point $p$ not on the hyperplane 
containing $Q$, 
and $C$ is the strict transform of a line passing through 
$p$ and disjoint from $Q$,
\item $Y$ is the blow-up of 
$\PP^{1}\times\PP^{n-1}$ at a point $p$
and $C$ is the strict transform of the fiber 
of the projection $\PP^{1}\times \PP^{n-1}\to 
\PP^{n-1}$
passing through $p$,
\item $Y$ is the blow-up of $\PP^{n}$ 
whose center is two distinct points $p$ and $q$, and
$C$ is the strict transform of the line $\overline{pq}$,
\item $Y$ is the blow-up of $\PP^{n}$ whose center is the 
union of two points $p,q$ and 
$\PP^{n-2}$ disjoint from $\overline{pq}$, and $C$ is the strict transform of $\overline{pq}$,
\item $Y$ is the blow-up of $\PP^{n}$ whose center is the 
union of two points $p,q$ and a smooth quadric 
$Q_{n-2}$ disjoint from $\overline{pq}$, and $C$ is the strict transform of $\overline{pq}$. 
\end{enumerate}
\end{theo}

\noindent {\it Remark.} We do not assume the projectivity of $Y$ because it follows from the assumption 
(see Lemma \ref{proj} below).

\section{Preliminaries}
We prove lemmas which will be needed for the proof of 
Theorem \ref{Fe=0}.

\begin{lem}\label{Ef=1}
We have $E\cdot f=1$.
\end{lem}
\begin{proof} 
Since $F\cdot e=0$, (the reduced part of) 
the intersection $E\cap F$ 
is a union of fibers of $\pi |_{E}:E\to C$. 
Hence $E\cap F$ is the exceptional locus  
of $\pi |_{F}:F\to \pi(F)$. 
Since $\varphi |_{F}:F\to W$ is a $\PP^{1}$-bundle, 
we see that  $E\cap F$ is a section of 
$\varphi |_{F}$. 
Hence we can write $E |_{F}=mE_{c}$ 
where $m$ is a natural number and 
$E_{c}:=\pi^{-1}(c)$ with $c\in C$ 
is a fiber of $\pi |_{E}$. 
Let $e_{c}$ be a line in $E_{c}\simeq \PP^{n-2}$. 
We have 
$$
mE_{c}\cdot e_{c}=
E |_{F}\cdot e_{c}=E\cdot e_{c}=-1, 
$$
where the first and second intersection numbers are
taken in $F$ and 
the last one is in $X$. 
Note that $(E_{c}\cdot e_{c})$ is an integer because 
$F$ is smooth. Thus we get $m=1$. 
It follows that $E \cap F$ is a reduced section of 
$\varphi |_{F}:F \to W$. Therefore, we have $E\cdot f=1$. 
\end{proof}

Now we consider $F_{Y}:=\pi(F)\subset Y$. Note 
that we have $\pi^{*}F_{Y}=F$.
\begin{lem}\label{fc=1}
We have $F_{Y}\cdot C=1$.
\end{lem}
\begin{proof}
Let $\widetilde{C}$ be a section of 
$\pi |_{E}:E\to C$. By (the proof of) 
 Lemma \ref{Ef=1}, $F|_{E}$ is a reduced 
 fiber of $\pi |_{E}$. Thus we have  
$F_{Y}\cdot C
=F_{Y}\cdot \pi_{*}\widetilde{C}
=F\cdot \widetilde{C}=1$. 
\end{proof}
By the proof of Lemma \ref{Ef=1}, 
we see that the intersection number 
$E_{c}\cdot e_{c}$ (taken in $F$) is equal to $-1$. 
It follows that $\pi |_{F} : F\to F_{Y}$
is the blow-up at the point $c$ whose 
exceptional divisor is $E_{c}$ and $F_{Y}$ is smooth. 
By the lemma \ref{Ef=1}, we see that $W$ is isomorphic to 
$E_{c}\simeq\PP^{n-2}$. We have the diagram:
\begin{equation*}
\begin{CD}
F @>\varphi |_{F} >> W\simeq \PP^{n-2} \\
@V\pi |_{F} VV \\
F_{Y} 
\end{CD}
\end{equation*}
where $\varphi |_{F}$ is a $\PP^{1}$-bundle. 
Note that $F$ is a Fano manifold. Indeed, we have 
$\rho(F)=2$ and $F$ has two extremal contractions 
$\pi |_{F}$ and $\varphi |_{F}$.   
According to the classification result from \cite{BCW}, 
this implies that $F_{Y}$ is isomorphic to $\PP^{n-1}$.
Furthermore, $f_{Y}:=\pi_{*}f$ is a line passing through 
the point $c=F_{Y}\cap C$. 
Note that $F_{Y}\cdot f_{Y}=F\cdot f=-1$. 
Hence there exists a blow-down 
$\varphi' :Y \to Y'$ contracting 
$F_{Y}\simeq \PP^{n-1}$ to a smooth point $p$, $Y'$
being (a priori) a complex manifold. 
Hence we have the commutative diagram:
\begin{equation}\label{diagram}
\begin{CD}
X @>\varphi >> Z \\
@V\pi VV @VV\pi ' V\\
Y @>>\varphi' > Y'
\end{CD}
\end{equation}
where $\pi' : Z \to Y'$ is the blow-up along the curve $C':=\varphi'(C)$. 

\begin{lem}\label{proj}
$Y$ is projective (hence, so is $Y'$).
\end{lem}
\begin{proof}
Assume to the contrary that $Y$ is not projective. 
Then, by \cite{BT} the normal bundle 
$N_{C/Y}$ is isomorphic to $\OO_{\PP^{1}}(-1)^{\oplus (n-1)}$. 
Since $\varphi' :Y\to Y'$ is a blow-up, $Y'$ is not projective either. 
Note that $Z$ is projective by the assumption of 
Theorem \ref{Fe=0}. 
Hence $N_{C'/Y'}\simeq \OO_{\PP^{1}}(-1)^{\oplus (n-1)}$ 
by \cite{BT} again. On the other hand, we have 
$N_{C/Y}\neq N_{C'/Y'}$ because $F_{Y}\cdot C>0$. 
Hence we get a contradiction. 
\end{proof}

We recall here the classification result due to 
\cite{BCW} which is indispensable to the proof of our 
Theorem \ref{Fe=0}.
Let $V_{d}$ denote the blow-up of $\PP^{n}$ 
along a smooth complete intersection 
$U_{d}:=H\cap D$ where $H$ is a hyperplane 
and $D$ is a hypersurface of degree $d$.
\begin{theo}[Bonavero, Campana and Wi\'sniewski \cite{BCW} Theorem 1.1]\label{bcw}
Let $Y'$ be a complex manifold of dimension $n\geq 3$. 
Let $\varphi': Y\to Y'$ be the blow-up at a point $p\in Y'$. 
Then $Y$ is a Fano manifold if and only if $Y'$ is 
isomorphic to either 
$\PP^{n}$, $Q_{n}$, or 
$V_{d}$ with $1\leq d\leq n$ 
and $p$ is not on the hyperplane $H$ containing the 
center $U_{d}$.
\end{theo}

\section{Proof of Theorem \ref{Fe=0}}
The proof is divided into two parts: 
\begin{itemize}
\item[(A)] If $Y$ is a Fano manifold, then $(Y,C)$ is one of the examples from (1) to (5).
\item[(B)] If $Y$ is not a Fano manifold, then $(Y,C)$ is one of the examples from (6) to (9).
\end{itemize}

Throughout the section, we frequently use the following: 
\begin{lem}[cf.\cite{MM} Proposition 7]\label{int}
Let $Y$ be a smooth projective variety of 
dimension $n\geq 3$, $C\subset Y$ 
a smooth subvariety of 
codimension $k\geq 2$, and $X$ the blow-up of $Y$ along $C$. 
Assume that $X$ is a Fano manifold. If 
$\Gamma\subset Y$ is a curve not contained in $C$ 
and $\Gamma\cap C\neq\varnothing$, then we have 
$(-K_{Y})\cdot \Gamma\geq k$.
\end{lem}
\begin{proof}
Let $\widetilde{\Gamma}$ be the strict transform 
of $\Gamma$ by the blow-up. For the exceptional 
divisor $E$, 
we have $E\cdot \widetilde{\Gamma}\geq 1$. 
Hence we have
$$
0<-K_{X}\cdot\widetilde{\Gamma}
=-K_{Y}\cdot \Gamma-(k-1)E\cdot\widetilde{\Gamma}
\leq -K_{Y}\cdot \Gamma-(k-1),
$$
which gives the statement.
\end{proof}

In what follows, we use the notation of the diagram (\ref{diagram}) in the previous section.
\subsection{Proof of (A)}
We assume that $Y$ is a 
Fano manifold. Since $\varphi':Y\to Y'$ is 
a blow-up at a point, we are exactly in the situation
of Theorem \ref{bcw}. 
Consider the extremal contraction 
$\gamma:Y\to Y''$ of ray $\R^{+}[g]$ such that
$F_{Y}\cdot g>0$ (see \cite{BCW} Lemme 2.1 for the existence of such a contraction). 
Then, by \cite{BCW} Proposition 2.2, $\gamma$ is either:
\begin{itemize}
\item[(A1)] a $\PP^{1}$-bundle, or 
\item[(A2)] a blow-up of a smooth projective variety along a smooth subvariety of codimension 2. 
\end{itemize}

In the case (A1),  
$Y'$ is isomorphic to $\PP^{n}$. We shall determine 
the position of $C$ in $Y$. If $C$ is not a fiber of 
$\gamma$, then there exists 
a fiber $\Gamma\simeq \PP^{1}$ of $\gamma$ 
such that $\Gamma\cap C\neq\varnothing$. 
Note that $-K_{Y}\cdot \Gamma=2$. Hence, 
by Lemma \ref{int} this is a contradiction. 
It follows that $C$ is a fiber of $\gamma$, i.e. 
the strict transform of a line in $Y'\simeq \PP^{n}$ 
passing through $p$, the center 
of the blow-up $\varphi'$. 
So, we get the example (1).

Now we treat the case (A2). Let $W_{\gamma}$ be 
the center of the blow-up $\gamma:Y\to Y''$ and $G$ the exceptional divisor. 
Note that $\gamma |_{G}:G\to W_{\gamma}$ 
is a $\PP^{1}$-bundle. 
By \cite{BCW}, there are two possibilities:
\begin{itemize}
\item[(A2]\hspace{-2mm}--1) $Y''$ is isomorphic to $\PP^{n}$ and $W_{\gamma}$ is isomorphic to $Q_{n-2}$, or 
\item[(A2]\hspace{-2mm}--2) $Y''$ is isomorphic to 
the $\PP^{1}$-bundle 
$\PP(\OO_{\PP^{n-1}}\oplus \OO_{\PP^{n-1}}(d-1))$ 
and $W_{\gamma}$ is a hypersurface in the section  
$\PP(\OO_{\PP^{n-1}})\simeq \PP^{n-1}$ 
whose normal bundle is isomorphic to 
$\OO_{\PP^{n-1}}(d-1)$. 
\end{itemize}

In the case (A2--1), 
$Y'$ is isomorphic to $Q_{n}$. 
Note that 
$F_{Y}=\Exc(\varphi')$ 
is the strict 
transform of the hyperplane containing 
$W\simeq Q_{n-2}$ 
by the blow-up $\gamma:Y\to Y''\simeq \PP^{n}$.
Since $F_{Y}\cdot C=1$ (Lemma \ref{fc=1}), $C$ is either 
a fiber of the $\PP^1$-bundle 
$\gamma |_{G}:G\to W$,  or 
the strict transform of a conic 
passing through $p\in Q_{n}$. 
So, we get the examples $(2)$ or (3).

In the case (A2--2), $Y'$ is isomorphic to $V_{d}$
with $1\leq d\leq n$. 
Let $\beta: Y'\simeq V_{d}\to \PP^{n}$ denote the blow-up along the smooth complete intersection 
$U_{d}=H\cap D$ with $H\in |\OO_{\PP^{n}}(1)|$ 
and $D\in |\OO_{\PP^{n}}(d)|$. 
Consider the composite of the two blow-ups
$\varepsilon:=\varphi'\circ \beta:Y\to \PP^{n}$. 
Note that the exceptional divisor $G$ of $\gamma$ 
is the strict transform by $\varepsilon$ 
of the cone over $U_{d}$ with vertex $\beta(p)$ 
(recall that $p$ is the center of the blow-up 
$\varphi':Y\to Y'$).  
Let $H_{Y}$ be the strict transform by $\varepsilon$ of 
the hyperplane $H$ containing $U_{d}$. 
We have $H_{Y}\cap F_{Y}=\varnothing$ because 
$\varepsilon(F_{Y})\notin H$ (see the statement of Theorem \ref{bcw}). 
\begin{claim}
We have $H_{Y}\cdot C=1$.
\end{claim}
\begin{proof}
Let $M$ be the exceptional divisor of the blow-up $\beta$ and 
$M_{Y}$ its strict transform by $\varphi'$. 
Since $F_{Y}\cdot C=1$ and 
$M_{Y}\cap F_{Y}=\varnothing$, 
we see that $C\not\subset M_{Y}$.  
In particular we have $\varepsilon_{*}C\not\equiv 0$. 
If $C\cap M_{Y}\neq \varnothing$, then there exists 
a fiber $\Gamma$ of the $\PP^{1}$-bundle 
$M_{Y}\to U_{d}$ meeting $C$. Note that 
$-K_{Y}\cdot \Gamma=1$.  By Lemma \ref{int}, this is a contradiction. Hence $M_{Y}\cdot C=0$. 
Note that $\varepsilon^{*}H=H_{Y}+M_{Y}$. 
We have 
$$
H_{Y}\cdot C=(H_{Y}+M_{Y})\cdot C
=\varepsilon^{*} H\cdot C
=H\cdot \varepsilon_{*}C>0.
$$ 
If $H_{Y}\cdot C\geq 2$ then there exists a line 
$h\subset H_{Y}\simeq \PP^{n-1}$ whose strict 
transform $\widetilde{h}$ by the blow-up 
$\pi:X\to Y$ satisfies 
$E\cdot \widetilde{h}\geq 2$. Then we have 
$$
K_{X}\cdot \widetilde{h}
=K_{Y}\cdot h+(n-2)E\cdot \widetilde{h}
\geq -n+d-1+2(n-2)=n+d-5\geq 0,
$$
which is a contradiction because $X$ is a 
Fano manifold. Hence we are done. 
\end{proof}
\begin{claim}
We have $d=1$ or $2$.
\end{claim}
\begin{proof}
Let $h$ be a line in $H_{Y}\simeq \PP^{n-1}$ 
such that $E\cdot \widetilde{h}=1$. Then we have 
$$
K_{X}\cdot \widetilde{h}
=K_{Y}\cdot h+(n-2)E\cdot \widetilde{h}
=d-3.
$$ 
Since $K_{X}\cdot \widetilde{h}<0$, we get $d=1$ or $2$.
\end{proof}
If $d=1$, we get the example (4) and if $d=2$, 
the example (5). The curve $C$ is determined by 
the condition $H_{Y}\cdot C=1$ and $F_{Y}\cdot C=1$.

\

\subsection{Proof of (B)}  Assume that $Y$ is not a 
Fano manifold. By \cite{Wis} Proposition 3.5, 
$E$ is isomorphic to $\PP^{1}\times \PP^{n-2}$ 
and $E\cdot l=-1$ where $l$ is a fiber 
of the projection 
$E\simeq \PP^{1}\times \PP^{n-2}\to \PP^{n-2}$.
Put $E_{Z}:=\varphi(E)$, $e_{Z}:=\varphi_{*}e$ and $l_{Z}:=\varphi_{*}l$. Since $E\cdot f=1$, 
$\varphi |_{E}:E\to E_{Z}$ is an isomorphism. 
Since $\varphi^{*}E_{Z}=E+F$, 
we have $E_{Z}\cdot e_{Z}=-1$ and 
$E_{Z}\cdot l_{Z}=0$.
Recall that $Y'$ is projective by Lemma \ref{proj}.
\begin{lem}
The projective varieties $Y'$ and $Z$ are  Fano manifolds. 
\end{lem}
\begin{proof} 
Since $E_{Z}\cdot l_{Z}=0$, we have 
$N_{C'/Y'}\simeq \OO_{\PP^{1}}^{\oplus (n-1)}\not\simeq  
\OO_{\PP^{1}}(-1)^{\oplus(n-1)}$. Therefore,  
$Y'$ is a Fano manifold by \cite{Wis} Proposition 3.5. 
Since $K_{X}=\varphi^{*}K_{Z}+F$, we get
$$
(-K_{Z})\cdot e_{Z}=-K_{X}\cdot e+F\cdot e=-K_{X}\cdot e>0.
$$
Note that the center $W$ of the blow-up 
$\varphi:X\to Z$ is a fiber of the projection 
$E_{Z}\simeq \PP^{1}\times\PP^{n-2}\to \PP^{1}$. 
Hence, any curve contained in $W$ is numerically proportional to a positive multiple of the line 
$e_{Z}$. By  
\cite{T2} Proposition 1, we conclude that 
$Z$ is a Fano manifold. 
\end{proof}

Since $Z$ is a Fano manifold, there exists an extremal ray 
$\R^{+}[m]\subset \overline{\NE}(Z)$ such that 
$E_{Z}\cdot m>0$ (see \cite{BCW} Lemme 2.1). We investigate 
the associated extremal contraction 
$\mu: Z\to Z'$.

\begin{lem}\label{curve}
We have $\mu_{*}B\not\equiv 0$ 
for any curve $B$ contained in $E_{Z}$.
\end{lem}
\begin{proof} Assume to the contrary that there exists a curve $B\subset E_{Z}$ such that $\mu_{*}B\equiv 0$. 
Then, there exists $a>0$ such that 
$B\equiv am$. On the other hand, 
we can write $B\equiv bl_{Z}+ce_{Z}$ 
with $b,c\geq 0$
because $B$ is contained in $E_{Z}\simeq \PP^{1}\times \PP^{n-2}$.
So, we have $E_{Z}\cdot B=E_{Z}\cdot(am)>0$ and 
$E_{Z}\cdot B=E_{Z}\cdot(bl_{Z}+ce_{Z})=-c\leq 0$,
a contradiction. \end{proof}
If there exists $z'\in Z'$ such that 
$\dim \mu^{-1}(z')\geq 2$, then there exists 
a curve $B$ contained in $E_{Z}\cap\mu^{-1}(z')$. 
Hence, any non-trivial fiber of $\mu$ 
has dimension at most $1$.
By \cite{A} (see \cite{Wis} Theorem 1.2),  
the extremal contraction $\mu$ is either: 
\begin{itemize}
\item[(B1)]  a conic bundle, or
\item[(B2)] a blow-up of a smooth projective variety along a smooth subvariety of codimension 2.
\end{itemize}

First we treat the case (B1).  We show that 
$\mu$ has no singular fiber, i.e. $\mu$ is a 
$\PP^1$-bundle.  Let $\Gamma$ be a 
fiber of $\mu$. Note that $\Gamma$ is 
isomorphic to $\PP^{1}$ or 
$\Gamma\simeq \Gamma_{1}\cup\Gamma_{2}$ 
with $\Gamma_{i}\simeq \PP^{1}\ (i=1,2)$. 

\begin{claim}\label{mult=1}
If $\Gamma$ meets $W$, then $\Gamma$ is a smooth fiber and the intersection $\Gamma\cap W$ is one point with multiplicity one.
\end{claim}
\begin{proof}
Assume  
$\Gamma=\Gamma_{1}\cup\Gamma_{2}$ 
with $\Gamma_{1}\cap W\neq\varnothing$. 
Note that $-K_{Z}\cdot \Gamma_{1}=1$. 
By Lemma \ref{int}, this is a contradiction. 
Hence $\Gamma$ is a smooth fiber.
Let $\widetilde{\Gamma}$ be the strict transform 
by $\varphi$. We have
$$
0<-K_{X}\cdot \wtilde{\Gamma}
=-K_{Z}\cdot \Gamma-F\cdot\widetilde{\Gamma}
= 2-F\cdot \widetilde{\Gamma},
$$
which gives $F\cdot \widetilde{\Gamma}=1$ and completes 
the proof.
\end{proof}
We conclude that $\mu |_{W}:W\to \mu(W)$ is an isomorphism. In particular, 
$\mu(W)\simeq\PP^{n-2}$. We put  
$M:=\mu^{-1}(\mu (W))$. Remark that  
$\mu |_{M}:M\to\mu (W)$ is a $\PP^{1}$-bundle and 
$W$ is a section. 
\begin{claim}
We have $E_{Z}\cdot \Gamma =1$.
\end{claim}
\begin{proof}
By Claim \ref{mult=1}, 
it is sufficient to prove $E_{Z}\cap M=W$. Let $\Gamma$ be any fiber of $\mu |_{M}:M\to \mu(W)$. We show that $E_{Z}\cap \Gamma\subset W$. 
Assume to the contrary that there exists a point 
$z\in E_{Z}\cap \Gamma$ such that $z\notin W$. 
Let $\Phi$ be the fiber of the $\PP^{n-2}$-bundle 
$\pi' |_{E_{Z}}:E_{Z}\to C'$ containing the point $z$. 
Since $\dim M\cap \Phi =n-3\geq 1$, there exists a curve $A\subset  M\cap \Phi$. 
Consider the ruled surface $S:=\mu^{-1}(\mu(A))$.
By Lemma \ref{curve} above, 
$\Gamma \not\subset E_{Z}$, 
hence $\pi' (\Gamma)\not\subset C'$ and we have $\dim\pi'(S)=2$. Therefore, 
$W\cap S$ and $A$ are exceptional curves on $S$.  
Note that $A\neq W\cap S$ because $\Phi\cap W=\varnothing$. Thus, we have a contradiction 
because $S$ is a ruled surface. 
\end{proof}

Now, we see that 
$\mu: Z \to Z'$ is a $\PP^{1}$-bundle and 
$\mu |_{E_{Z}}:E_{Z}\to Z'$ is an isomorphism. 
It follows that $Z'$ is isomorphic to 
$\PP^{1}\times\PP^{n-2}$. 
 Pushing down the exact sequence: 
 $$
 0\to \OO_{Z}\to \OO_{Z}(E_{Z}) \to \OO_{E_{Z}}(E_{Z})\to 0,
 $$
 we get
 $$
 0\to \mu_{*}\OO_{Z}\to \mu_{*}\OO_{Z}(E_{Z})\to \mu_{*}\OO_{E_{Z}}(E_{Z})\to 
 R^{1}\mu_{*}\OO_{Z}=0.
 $$
Since $\mu$ is an extremal contraction, we have $\mu_{*}\OO_{Z}\simeq \OO_{\PP^{1}\times \PP^{n-2}}$. 
Recall that 
$\OO_{E_{Z}}(E_{Z})\simeq 
\OO_{\PP^{1}\times \PP^{n-2}}(0,-1)$.  
Since $\mu |_{E_{Z}}$ is an isomorphism, 
we have $\mu_{*}\OO_{E_{Z}}(E_{Z})\simeq \OO_{\PP^{1}\times \PP^{n-2}}(0,-1)$. Thus we get 
the splitting sequence
$$
0\to \OO_{\PP^{1}\times \PP^{n-2}}\to \mu_{*}\OO_{Z}(E_{Z})\to 
\OO_{\PP^{1}\times\PP^{n-2}}(0,-1)\to 0, 
$$
which gives 
$\mu_{*}\OO_{Z}(E_{Z})\simeq \OO_{\PP^{1}\times \PP^{n-2}}\oplus \OO_{\PP^{1}\times\PP^{n-2}}(0,-1)$. Thus we have
$$
Z\simeq \PP (\OO_{\PP^{1}\times \PP^{n-2}}\oplus \OO_{\PP^{1}\times\PP^{n-2}}(0,-1))
\simeq \PP^{1}\times \Bl_{p}(\PP^{n-1}), 
$$
where $\Bl_{p}(\PP^{n-1})$ denotes the blow-up 
of $\PP^{n-1}$ at the point $p$.
We see that $Y'\simeq \PP^{1}\times \PP^{n-1}$ and $C'$ is a fiber of the 
projection $Y'\to \PP^{n-1}$. We obtain the example (6).

\

Now we consider the case (B2).
Let $F_{Z}$ be the exceptional divisor of 
the blow-up $\mu : Z \to Z'$. Since $E_{Z}$ 
is strictly positive on the extremal ray $\RR[m]$, 
we have $E_{Z}\neq F_{Z}$, in particular 
$F_{Z}\cdot e_{Z}\geq 0$.
If $F_{Z}\cdot e_{Z}>0$, there exists a fiber 
$m_{0}$ of the $\PP^{1}$-bundle 
$F_{Z}\to \mu(F_{Z})$ such that 
$m_{0}\cap W\neq \varnothing$ (recall that $W$ denote the center of the blow-up $\varphi : X \to Z$). Since 
$-K_{Z}\cdot m_{0}=1$, we get a contradiction 
by Lemma \ref{int}. Hence we have $F_{Z}\cdot e_{Z}=0$. 

Recall that 
$\pi':Z \to Y'$ is the blow-up along $C'$ and  
$\mu : Z \to Z'$ is a blow-up along a center of codimension 2 with $F_{Z}\cdot e_{Z}=0$.
Since $Y'$ and $Z$ are Fano manifolds, we can use the statement (A) (already proved in the previous subsection)
to classify the pairs $(Y',C')$. Moreover, we have the condition on the normal bundle: 
$N_{C'/Z'}\simeq \OO_{\PP^{1}}^{\oplus(n-1)}$, 
which is satisfied for the following 
cases: 
\begin{itemize}
\item $Y'$ is the blow-up of $\PP^{n}$ at a point $q$ and $C'$ is the strict transform of the line 
passing through $q$
\item $Y'$ is the blow-up of 
$\PP^{n}$ at a point $q$ and a linear subspace 
$P\simeq \PP^{n-2}$ and $C'$ is the strict transform 
of a line passing through $q$, 
\item $Y'$ is the blow-up of $\PP^{n}$ 
at a point $q$ and a quadric $Q\simeq Q_{n-2}$
 and $C'$ is the strict transform of 
a line passing through $q$.
\end{itemize}
Recall that in each case, $Y$ is the blow-up of $Y'$ at the point $p\in C'$. So, we get the examples (7), (8), and (9).
Hence, the proof of the statement (B) is completed.

\noindent -----------------------------------------

\noindent {\small Toru TSUKIOKA \ \ \ \ e-mail:\ tsukiokatoru@yahoo.co.jp\\ 
Faculty of Liberal arts and Sciences, 
Osaka Prefecture University \\
1-1\ Gakuen-cho\ Nakaku\ Sakai, 
Osaka 599-8531\ Japan
}

\end{document}